% Logic Eprints
%Submitted 1213 Thu Sep 09, 1993 by: brendle@bimacs.cs.biu.ac.il (joerg brendle)
%logic/brendle/poro.tex
%

\magnification=\magstep1

%Das Folgende sollte die "Blackboard bold font" \Bbb definieren, und 
%auch  die "Fraktur" font, \frak:
% Beispiel:    $ {\frak c} = |{\Bbb R}| $

\def\hexnumber#1{\ifcase#1 0\or1\or2\or3\or4\or5\or6\or7\or8\or9\or
	A\or B\or C\or D\or E\or F\fi }

%  The following lines establish the use of the Euler Fraktur font.
\font\teneuf=eufm10
\font\seveneuf=eufm7
\font\fiveeuf=eufm5
\newfam\euffam
\textfont\euffam=\teneuf
\scriptfont\euffam=\seveneuf
\scriptscriptfont\euffam=\fiveeuf

%  End definition of Euler Fraktur font.

\font\tenmsx=msam10
\font\sevenmsx=msam7
\font\fivemsx=msam5
\font\tenmsy=msbm10
\font\sevenmsy=msbm7
\font\fivemsy=msbm5
\newfam\msxfam
\newfam\msyfam
\textfont\msxfam=\tenmsx  \scriptfont\msxfam=\sevenmsx
  \scriptscriptfont\msxfam=\fivemsx
\textfont\msyfam=\tenmsy  \scriptfont\msyfam=\sevenmsy
  \scriptscriptfont\msyfam=\fivemsy
\edef\msx{\hexnumber\msxfam}

\mathchardef\upharpoonright="0\msx16
\let\restriction=\upharpoonright
\def\Bbb#1{\tenmsy\fam\msyfam#1}

\def\restrict{{\restriction}}

\def\Smallskip{\vskip1.4truecm}
\def\Bigskip{\vskip2.2truecm}

\def\qed{{\vcenter{\hrule height.4pt \hbox{\vrule width.4pt height5pt
 \kern5pt \vrule width.4pt} \hrule height.4pt}}}
\def\notin{{\in}\kern-5.5pt / \kern1pt}

\def\II{{\Bbb I}}

\def\sm{{\smallskip}}
\def\ce#1{{\centerline{#1}}}
\def\no{{\noindent}}
\def\la{{\langle}}
\def\ra{{\rangle}}
\def\sub{\subseteq}

\font\small=cmr8 scaled\magstep0
\font\smalli=cmti8 scaled\magstep0
\font\capit=cmcsc10 scaled\magstep0
\font\capitg=cmcsc10 scaled\magstep1

\font\sanse=cmss10 scaled\magstep0

\font\bolds=cmssdc10 scaled\magstep0

\overfullrule=0pt
\openup1.5\jot

\centerline{\capitg The additivity of porosity ideals}
\footnote{}{{\openup-6pt {\small {\smalli 1991 
Mathematics subject classification.}
03E05 \par {\smalli Key words and
phrases.} porous sets, cardinal invariants
\endgraf}}}
\Bigskip
\centerline{\sanse J\"org Brendle\footnote{$^*$}{{\small
The author wishes to thank the MINERVA-foundation
for supporting him}}}
\Smallskip
{\baselineskip=0pt {\small
\noindent 
Department of Mathematics,
Bar--Ilan University,
52900 Ramat--Gan, Israel
\par
and \par
\noindent Mathematisches Institut der Universit\"at T\"ubingen,
Auf der Morgenstelle 10, 72076 T\"ubingen, Germany}} 
\Bigskip
{\narrower 
\no{\sanse Abstract.} We show that several $\sigma$--ideals related
to porous sets have additivity $\omega_1$ and cofinality $2^\omega$.
This answers a question addressed by Miroslav Repick\'y. \par}
\Bigskip
\no {\bolds 1. Introduction.} Given a $\sigma$--ideal ${\cal I}$
on the reals $\II := [0,1]$ we let
\sm
$add({\cal I}) := $ the least $\kappa$ so that $\exists 
{\cal F} \in [{\cal I}]^\kappa \; (\cup {\cal F} \not\in
{\cal I})$;
\par
$cov({\cal I}) :=$ the least $\kappa$ so that $\exists
{\cal F} \in [{\cal I}]^\kappa \; (\cup {\cal F} = \II$);
\par
$unif({\cal I}) := $ the least $\kappa$ so that $[\II]^\kappa
\setminus {\cal I} \neq \emptyset$;
\par
$cof({\cal I}) :=$ the least $\kappa$ so that $\exists {\cal F}
\in [{\cal I}]^\kappa \; \forall A \in {\cal I} \; \exists
B \in {\cal F} \; (A \sub B)$.
\sm
\no Next, let $A\sub \II$ be a set of reals. The {\it porosity}
and the {\it symmetric porosity} of the set $A$ at a real
$r\in\II$ are defined by
$$p(A,r) = \lim\sup_{\epsilon \to 0^+} {\lambda (A, (r - 
\epsilon , r + \epsilon)) \over \epsilon}\; , \;{\rm and}
$$
$$s(A,r) = \lim\sup_{\epsilon \to 0^+} {\lambda^* (A, (r -
\epsilon , r + \epsilon )) \over \epsilon} \; , $$
respectively, where $\lambda (A,I)$ denotes the maximal
length of an open subinterval of the interval $I$
which is disjoint from $A$; similarly $\lambda^* (A,(c,d))$
stands for the maximal $\delta \geq 0$ so that $(c, c+ \delta)
\cup (d-\delta , d)$ is disjoint from $A$.
For $a\in A$ we have $s(A,a) \leq p(A,a) \leq 1$.
$A$ is {\it porous} ({\it strongly porous}, respectively)
iff $p(A,a) > 0$ ($p(A,a) = 1$, resp.) for every $a\in A$;
similarly, $A$ is {\it symmetrically porous} ({\it strongly
symmetrically porous}, respectively) iff $s(A,a) > 0$
($s(A,a) = 1$, resp.) for every $a\in A$. Finally, we
let 
\sm
${\cal P}:=$ the $\sigma$--ideal generated by the strongly porous
sets;
\par
${\cal P}^+ :=$ the $\sigma$--ideal generated by the porous sets;
\par
${\cal S} :=$ the $\sigma$--ideal generated by the strongly
symmetrically porous sets;
\par
${\cal S}^+ := $ the $\sigma$--ideal generated by the symmetrically
porous sets.
\sm\par
\no The elements of ${\cal P}$ (${\cal P}^+$, ${\cal S}$, ${\cal
S}^+$, respectively) are called {\it $\sigma$--strongly porous
sets} ({\it $\sigma$--porous sets}, etc., respectively).
\bigskip
The goal of this note is to show (in section 2) $add ({\cal P})
= add ({\cal P}^+) = add ({\cal S}) = add ({\cal S}^+) =
\omega_1$ and $cof ({\cal P}) = cof ({\cal P}^+) = cof ({\cal S})
= cof ({\cal S}^+) = 2^\omega$. In section 3 we make some
remarks, and state a question, concerning the invariants
$unif ({\cal I}) $ and $cov ({\cal I})$ where ${\cal I}$
is again one of the four ideals defined above.
For a survey concerning known results about these
cardinals, as well as further references, we refer the
reader to the survey article [Re].
\bigskip
\bigskip
\no {\bolds 2. Proof of the main result.} We shall work in a
space of the form $X := \prod_n g(n)$, where $g \in (\omega
\setminus 2)^\omega$; i.e. $f \in X$ iff $f\in\omega^\omega$
and for all $n \in\omega$, $f(n) < g(n)$. Using the Cantor
expansion, we can almost identify $X$ and the unit interval $\II$
in a canonical way: reals $r \in\II$ correspond to
reals $f\in X$ via the map $\phi : X \to \II$ defined
by
$$r = \phi (f) = \sum_{n\in\omega} {f(n) \over g(0) \cdot
... \cdot g(n)}\;\; .$$
However, we shall be more interested in a slightly different
identification. Let $g\in\omega^\omega$ be strictly increasing
and taking odd values such that there are strictly increasing
sequences of natural numbers $\la m_n ;\; n\in\omega\ra$ and
$\la \ell_n ;\; n\in\omega\ra$ so that
\sm
\item{(i)} ${\ell_n^{m_n} \over g(n) } \to 0$ for $n \to
\infty$;
\par
\item{(ii)} ${\ell_{n-1} \over \ell_n} \to 0$ for $n \to \infty$.
\sm
\no Next, given $n\in\omega$, we let $\mu_n$ be the unique
measure on $P(g(n))$ satisfying
\sm
\item{(a)} $\mu_n (g(n)) = 1$;
\par
\item{(b)} $\mu_n (\{ i \}) = \mu_n (\{ j \}) $ for
$m_n \leq i \leq j \leq g(n) - m_n -1$;
\par
\item{(c)} $\mu_n (\{ i \}) = \mu_n (\{ g(n) - i -1\})
= {\mu_n (\{ i+1 \} ) \over \ell_n }$ for $i < m_n$.
\par
\sm
\no We equip $X$ with the product measure of the $\mu_n$.
This gives another almost--identification of $X$ and $\II$:
the idea is that for $\sigma \in \prod_{n<\ell} g(n)$,
the open set $[\sigma] := \{ f \in X ; \; \sigma \sub
f \}$ corresponds to an interval in $\II$ of length
$\prod_{n<\ell} \mu_n (\{ \sigma (n) \})$. We call this
almost--correspondence $\tilde\phi : X \to \II$, and
say $A\sub X$ is porous (or strongly porous, etc.)
iff $\tilde \phi (A)$ is. Since $\tilde\phi \circ \tilde\phi^{-1}
(B) = B$ and $\tilde\phi^{-1} \circ \tilde\phi
(A) = A \cup C$ for some at most countable $C\sub X$, the
$\sigma$--ideals of $\sigma$--porous ($\sigma$--strongly porous, etc.)
sets can be identified.
\sm
Choose $\la k_n ;\; n\in\omega \ra$ a strictly increasing
sequence of natural numbers so that $k_{n+1} - k_n
\geq 10$. Let $\la A_\alpha ; \; \alpha < 2^\omega \ra$
be a sequence of almost--disjoint subsets of $\omega$.
For $\alpha < 2^\omega$ we define $B_\alpha : = \{
f\in X ;\; \forall n \in A_\alpha \; (f(k_n) = {g(k_n) - 1 \over
2} ) \}$. Each $B_\alpha$ is easily seen to be strongly
symmetrically porous. We claim:
\medskip
{\capit Theorem.} {\it Whenever $C \sub X$ is $\sigma$--porous
then for all but countably many $\alpha \in 2^\omega$,
$B_\alpha \not\sub C$.}
\medskip
{\capit Corollary.} {\it For ${\cal I} \in \{ {\cal P} ,
{\cal P}^+ , {\cal S} , {\cal S}^+ \}$, $add ({\cal I})
=\omega_1$ and $cof ({\cal I}) = 2^\omega$.}
\medskip
{\it Proof of Corollary from Theorem.} As ${\cal S}
\sub {\cal S}^+ , {\cal P} \sub {\cal P}^+$, $\{ B_\alpha ;
\; \alpha \in \omega_1 \}$ witnesses $add ({\cal I}) =
\omega_1$ for any of the $\sigma$--ideals. Furthermore any
family ${\cal F}$ of sets from ${\cal I}$ of size $< 2^\omega$
cannot be cofinal, because some $B_\alpha$ will not be contained in
any member of ${\cal F}$. $\qed$
\medskip
{\it Proof of Theorem.} Let $C = \bigcup_{i\in\omega} C_i$,
where each $C_i$ is porous. Fix $\sigma \in \prod_{n<\ell}
g(n)$ (for some $\ell \in\omega$), $m < [ {\ell \over 2 } ]$
and $\Gamma \sub 2^\omega$ finite. We set
\sm
\ce{$B (\sigma , m , \Gamma ) := \{ f\in X ; \; \sigma \sub f \;\land\;
\forall n \geq \ell \; (m\leq f(n) \leq g(n) - m-1 ) \;\land\; $}
\par
\ce{$\;\land\;
\forall \alpha\in\Gamma
\;\forall n \in A_\alpha \; (k_n \geq \ell \to f(k_n) = {g(k_n)
-1\over
2} ) \}$ and}
\medskip
\ce{$\hat B (\sigma ,m , \Gamma ) := \{ \tau \supseteq 
\sigma ;\; \forall n \; (\ell \leq n < lh (\tau) \to m \leq
\tau (n) \leq g(n) - m-1 ) \;\land\; $}
\par
\ce{$\;\land\;\forall \alpha\in\Gamma\;\forall
n \in A_\alpha \; (\ell \leq k_n < lh (\tau) \to \tau (k_n) =
{g(k_n) - 1 \over 2 } ) \}$.}
\medskip
\no Given $i\in\omega$, we say $\Gamma$ is {\it $(\sigma , m
,i)$--funny} iff 
\sm
\itemitem{{\sanse either (I)}} there are uncountably many $\Delta_\alpha$
($\alpha < \omega_1$) which are pairwise disjoint
so that for all $\alpha < \omega_1$
\par
\ce{$B(\sigma , m+ 1 , \Gamma \cup \Delta_\alpha ) \cap C_i
= \emptyset$;}
\sm
\itemitem{{\sanse or (II)}} for some $\tau \supseteq \sigma \;
([\tau] \cap C_i = \emptyset \;\land\; \tau \in \hat B (\sigma
, m , \Gamma ) )$.
\medskip
{\capit Main Observation.} {\it Given $\sigma , m, i$ as above,
there is $\Omega \sub 2^\omega$ countable so that whenever
$\Delta \sub 2^\omega \setminus \Omega$ is finite, then
$\Delta$ is $(\sigma , m ,i)$--funny.}
\medskip
{\it Proof of Main Observation.}
Suppose not. Then we can easily construct a sequence $\la
\Delta_\alpha ;\; \alpha < \omega_1 \ra$ of pairwise disjoint finite
sets none of which is $(\sigma , m , i)$--funny. By clause (I)
of the definition of funniness applied to $\Delta_0$ there
is $\alpha < \omega_1$ so that $B (\sigma , m+1 , \Delta_0 \cup
\Delta_\alpha ) \cap C_i \neq \emptyset$. Choose $f$ from the
latter set. Let $p := p (\tilde\phi (C) , \tilde\phi (f)) 
>0$. Find $\epsilon$ so small that $\epsilon$ and $n = n_\epsilon$,
which is unique with $\mu ( [f\restrict n ]) \geq\epsilon$ and
$\mu ( [ f\restrict n+1 ]) < \epsilon$, satify:
\sm
\item{(A)} $\lambda := \lambda (\tilde\phi (C_i) ,
(\tilde\phi (f) - \epsilon , \tilde\phi (f) + \epsilon
)) > {p \over 2 } \cdot \epsilon$;
\par
\item{(B)} ${\ell_n^{m_n} \over g(n)}$ is small compared to
$p$;
\par
\item{(C)} ${\ell_{n-1} \over \ell_n}$ is small compared to $p$;
\par
\item{(D)} $( \bigcup_{\beta \in \Delta_0} A_\beta )
\cap (\bigcup_{\beta \in \Delta_\alpha} A_\beta) \sub n'$
and $n \geq k_{n'}$.
\sm
\par
\no Without loss $I := ( \tilde\phi (f) - \epsilon , \tilde\phi
(f) - \epsilon + \lambda )$ is disjoint from $\tilde\phi (C_i)$.
Clearly $\tilde\phi (f) - \epsilon \not\in \tilde\phi 
[f\restrict n+1 ]$. Also we either have $\tilde\phi (f) - \epsilon
\in \tilde\phi [f\restrict n \hat{\;} \la j\ra ]$ for some $j <
f(n)$
or $\tilde\phi (f) -\epsilon \in\tilde\phi [f\restrict (n-1)
\hat{\;} \la j\ra ]$ for some $j < f(n-1)$ or $\tilde\phi (f)
-\epsilon \in\tilde\phi [f \restrict (n-2) \hat{\;} \la f(n-2)
-1 \ra ]$ [this is because ${\mu ([ f\restrict n-1 ])
\over g(n-1) - 2 \cdot m_{n-1}} > \mu ([ f\restrict n ])
\geq \epsilon$, which implies $\mu ([ f\restrict (n-2) \hat{\;}
\la f (n-2) -1 \ra ]) > \epsilon \cdot {g(n-1) - 2 \cdot
m_{n-1} \over \ell_{n-2} } > \epsilon$].
The core of the proof is to
construct $\tau \supseteq f \restrict n-2$ so that 
$lh (\tau) \leq n +2$, $\forall j \; (n-2 \leq j \leq n+1 
\to m \leq \tau (j) \leq g(j) - m -1)$ and $[\tau] \cap C_i
= \emptyset$ $(\clubsuit)$.
\par
It is easy to see that one of the following three cases
must hold:
\sm
\no{\sanse Case 1. For $\sigma_j = f \restrict j-1$, where
$j\in\{ n,n+1 \}$, we have $\tilde\phi (f) - \epsilon \in\tilde\phi
[\sigma_j \hat{\;} \la k \ra ]$, where $m\leq k < f(j-1)$
($k < f(j-1)-1$ in case $j=n$).}
\par
\item{} Note that for any such $k$ we have $\mu ([ \sigma_j \hat{\;}
\la k \ra ]) \leq \mu ([ \sigma_j \hat{\;} \la k+1 \ra ]) < \epsilon
\cdot \ell_{j-1}^{m_{j-1}}$ [in case $j=n$, this follows because the
assumption $\tilde\phi (f) - \epsilon \in \tilde\phi [\sigma_n
\hat{\;} \la k \ra ]$ implies $\mu ([ \sigma_n \hat{\;}
\la f (n-1) - 1 \ra ]) < \epsilon$; in case $j= n+1$, this is
immediate from $\mu ([ f\restrict n+1 ]) < \epsilon$].
Let $\ell$ be so that $\tilde\phi (f) - \epsilon \in\tilde\phi
[\sigma_j \hat{\;} \la k \ra \hat{\;} \la \ell \ra ]$. In case
$\ell < g(j) -m-1$ find $0< \ell ' \leq m$ so that $m \leq \ell
+ \ell ' \leq g(j) -m-1$, and we have
$$\sum_{\tilde\ell \leq \ell '} \mu ([ \sigma_j \hat{\;}
\la k \ra \hat{\;} \la \ell + \tilde\ell\ra ]) < \epsilon 
\cdot { (m+1) \cdot \ell_{j-1}^{m_{j-1}} \over g(j) - 2
\cdot m_j } < \epsilon \cdot {p\over 2} < \lambda.$$
This entails $\tilde\phi [\sigma_j \hat{\;}
\la k\ra\hat{\;} \la \ell + \ell ' \ra ] \sub I$, and
so $[\sigma_j \hat{\;}\la k\ra\hat{\;}\la \ell+ \ell '\ra] \cap
C_i = \emptyset$. Hence $\tau = \sigma_j \hat{\;} \la k\ra
\hat{\;} \la\ell + \ell '\ra$ will work. In case $\ell \geq g(j)
-m-1$ we compute
$$\sum_{g(j) > \tilde\ell \geq \ell} \mu ([ \sigma_j \hat{\;}
\la k \ra \hat{\;} \la \tilde\ell \ra ]) + \sum_{\tilde\ell
\leq m} \mu ([\sigma_j \hat{\;} \la k+1 \ra \hat{\;} \la
\tilde\ell\ra ]) < \epsilon \cdot {2 \cdot (m+1) \cdot
\ell_{j-1}^{m_{j-1}} \over g(j) - 2 \cdot m_j } < \epsilon
\cdot {p\over 2} < \lambda.$$
Thus $\tau = \sigma_j \hat{\;} \la k+1 \ra
\hat{\;} \la m\ra$ is as required.
\par
\sm
\no{\sanse Case 2. For $\sigma_j = f\restrict j-2$, where
$j \in\{ n,n+1\}$, we have:
\par
\no either $\tilde\phi (f) - \epsilon \in \tilde\phi
[\sigma_j \hat{\;} \la f (j-2) -1 \ra \hat{\;} \la g(j-1) -m-1
\ra \hat{\;} \la g(j) - k - 1 \ra ]$, where $k \leq m$,
\par
\no or $\tilde\phi (f) - \epsilon \in\tilde\phi [\sigma_j
\hat{\;} \la f(j-2) -1 \ra \hat{\;} \la g(j-1) - k-1 \ra ]$,
where $k < m$,
\par
\no or $\tilde\phi (f) - \epsilon \in\tilde\phi
[f\restrict (j-1) \hat{\;} \la k \ra ]$, where $k<m$.}
\par
\item{} In this case we necessarily have $\mu ([ f \restrict
(j-1) \hat{\;} \la m \ra ] ) < \epsilon$, and hence
$\mu ([ \sigma_j \hat{\;} \la f(j-2)-1\ra \hat{\;} \la g(j-1) -m-1
\ra ]) < \epsilon \cdot \ell_{j-2}$.
Thus (putting $\tilde\sigma_j = \sigma_j \hat{\;} \la f(j-2) - 1
\ra$ and $\hat\sigma_j = f \restrict j-1$)
$$\sum_{k\leq m} \mu ([\tilde\sigma_j \hat{\;}  \la
g(j-1) -m-1 \ra \hat{\;} \la g(j) - k -1 \ra ])
+ \sum_{k<m} \mu ([ \tilde\sigma_j 
\hat{\;} \la g(j-1)-k-1 \ra ]) +$$
$$+ \sum_{k<m} \mu ([ \hat\sigma_j \hat{\;} \la k \ra ])
+ \sum_{k\leq m} \mu ([\hat\sigma_j \hat{\;} \la m \ra \hat{\;}
\la k \ra ]) < \epsilon \cdot 2\cdot( {(m+1) \cdot \ell_{j-2} \over g(j)
- 2 \cdot m_j} +\sum_{k<m} {\ell_{j-2} \over \ell_{j-1}^{k+1}})
< \epsilon \cdot {p\over 2} <\lambda.$$
Hence $\tau = f \restrict (j-1) \hat{\;} \la m\ra \hat{\;} \la m
\ra$ is as required.
\sm
\par
\no{\sanse Case 3. For $\sigma_j = f\restrict (j-2) \hat{\;}
\la f(j-2)-1 \ra$, where $j\in \{ n,n+1 \}$, we have:
\par
\no either $\tilde\phi (f) -\epsilon \in \tilde\phi
[\sigma_j \hat{\;} \la k \ra ]$, where $k < g(j-1) - m -1$,
\par
\no or $\tilde\phi (f) - \epsilon \in\tilde\phi [\sigma_j
\hat{\;} \la g(j-1) - m-1 \ra \hat{\;} \la k \ra ]$, where 
$k < g(j) -m-1$.}
\par
\item{} As before we necessarily have $\mu ([ f\restrict (j-1)
\hat{\;} \la m \ra ]) < \epsilon$, and hence $\mu ([ \sigma_j
\hat{\;} \la k \ra ]) < \epsilon \cdot \ell_{j-2} \cdot
\ell_{j-1}^{m_{j-1}}$ (for $k \leq g(j-1) - m-1$). If $k \geq m$
we can finish similarly to case 1
[i.e. we let $\ell$ be so that $\tilde\phi (f)
-\epsilon \in \tilde\phi [\sigma_j \hat{\;} \la k \ra
\hat{\;} \la\ell\ra ]$ and split into the two subcases
$\ell < g(j) -m-1$ and $\ell \geq g(j) -m-1$]. Note that
$$\sum_{\ell < m} \mu ([ \sigma_j \hat{\;} \la \ell \ra ])
+\sum_{\tilde \ell \leq m} \mu ([ \sigma_j \hat{\;}
\la m\ra \hat{\;}\la \tilde\ell\ra ]) < \epsilon \cdot (\sum_{\ell <
m} {\ell_{j-2} \over \ell_{j-1}^{\ell + 1} } + {(m+1) \cdot
\ell_{j-2} \over g(j) - 2 \cdot m_j} ) <\epsilon\cdot{p\over 2}
<\lambda.$$
Thus if $k<m$, $\tau = \sigma_j \hat{\;} \la m\ra \hat{\;} \la m
\ra$ will work (in fact this is similar to case 2).
\sm\par
\no Thus we have found $\tau$ as required in $(\clubsuit)$.
By almost--disjointness and the choice of $n$ we necessarily
have either $\tau \in \hat B(\sigma , m , \Delta_0)$
or $\tau \in \hat B (\sigma , m , \Delta_\alpha )$.
Thus either $\Delta_0$ or $\Delta_\alpha$ is funny,
and we reach a contradiction. This proves the
Main Observation. $\qed$
\medskip
{\it Conclusion.} By the Main Observation, there is $\Omega
\sub 2^\omega$ countable so that whenever $\Gamma \sub
2^\omega \setminus \Omega$ is finite, then $\Gamma$ is
$(\sigma , m ,i)$--funny for all $\sigma \in \prod_{n<\ell}
g(n)$ (some $\ell\in\omega$), $m < [{\ell \over 2}]$
and $i\in\omega$.
\par
Fix $\alpha \in 2^\omega \setminus \Omega$. We construct recursively
$\sigma_i \in \prod_{n < \ell(i)} g(n)$, where $\ell (i) \in
\omega$, $m(i) < [{\ell (i) \over 2 }]$ and $\Gamma_i \sub 2^\omega
\setminus \Omega$ finite so that 
\sm
\item{(1)} $\{\alpha\} = \Gamma_0$, $\Gamma_i \sub \Gamma_{i+1}$;
\par
\item{(2)} $\ell (0) = 2$ and $\ell (i) < \ell (i+1)$;
\par
\item{(3)} $m(0) = 0$ and $m(i) \leq m (i+1)$;
\par
\item{(4)} $
\sigma_{i+1} \in \hat B (\sigma_i , m(i) , \Gamma_i ) $;
\par
\item{(5)} $B(\sigma_{i+1} , m (i+1) , \Gamma_{i+1} ) \cap
C_i = \emptyset$.
\par
\sm
\no --- $i=0$. Trivial.
\par
\no --- $i \to i+1$. $\Gamma_i$ is $(\sigma_i , m(i) , i)$--funny.
So either find $\Gamma_i \sub \Gamma_{i+1} \sub 2^\omega
\setminus \Omega$ so that $B (\sigma_i , m(i) + 1 , \Gamma_{i+1}
) \cap C_i = \emptyset$; in this case $\ell (i+1) = \ell (i) + 2$,
$\sigma_{i+1} \in \hat B (\sigma_i , m(i) + 1 , \Gamma_{i+1})$
of length $\ell (i+1)$, $m (i+1) = m(i) + 1$, and (1) --- (5)
are satisfied. Or find $\tau \supseteq \sigma_i$ so that 
$[\tau] \cap C_i = \emptyset$ and $\tau \in \hat B (\sigma_i
, m(i) , \Gamma_i )$. In case $\tau \supset \sigma_i$
let $\sigma_{i+1} = \tau$ and $\ell (i+1) = lh (\tau)$;
otherwise $\ell (i) = \ell (i) + 1$ and $\sigma_{i+1}
\in \hat B (\sigma_i , m(i) , \Gamma_i )$. We let 
$\Gamma_{i+1} = \Gamma_i$ and $m (i+1) = m(i)$.
Again (1) --- (5) are satisfied.
\par
This concludes the construction.
Let $f = \bigcup_{i\in\omega} \sigma_i$. It is easily seen that
$f\in B_\alpha \setminus C$, thus proving the Theorem.
$\qed$
\bigskip
\bigskip
\no {\bolds 3. Porosity and evasion.} We were motivated to
prove the above result by our discussion of evasion ideals
[Br, section 3] which seem to be closely related to porosity
ideals.
\par
As in section 2,
fix $g\in (\omega \setminus 2)^\omega$, and let $X := \prod_n
g(n)$. Following Blass [Bl, section 4] (see also [Br, 3.1.]),
an {\it $X$--predictor} is a pair $\pi = (D_\pi ,  (\pi_n ;\;
n \in D_\pi ))$ such that for every $n \in D_\pi$,
$\pi_n : \prod_{k<n} g(k) \to g(n)$; $\pi$  {\it predicts}
$f\in X$ iff $\forall^\infty n \in D_\pi \; (\pi_n (f
\restrict n )  = f(n))$; otherwise $f$ {\it evades}
$\pi$. ${\bf e}_X := min \{ \vert {\cal F} \vert ; 
\; {\cal F}\sub X \;\land\; \forall X$--predictors $\pi \;
\exists f\in X \; (f$ evades $\pi ) \}$ is the {\it evasion
number}.
Furthermore, let ${\cal I}_X := \{ A \sub X ;\; $ there is
a countable set of $X$--predictors $\Pi$ so that for all
$f\in A$ there is $\pi\in\Pi$ predicting $f \}$ [Br, 3.5.].
Making again a standard identification between $X$ and
$\II$ as at the beginning of section 2, we see that 
${\cal I}_{2^\omega} \sub {\cal P}^+$ and
${\cal I}_X \sub {\cal P}$ for $X = \prod_n g(n)$,
where $g$ converges to infinity [moreover, the sets
$B_\alpha$ which are crucial for the proof in section 2
are elements of ${\cal I}_X$].
Thus
${\bf e}_X \leq unif ({\cal I}_X) \leq unif ({\cal P})$ and
${\bf e}_{2^\omega} \leq unif ({\cal I}_{2^\omega}) \leq
unif ({\cal P}^+)$ as well as $cov({\cal P}) \leq cov ({\cal I}_X)$
and $cov({\cal P}^+) \leq cov ({\cal I}_{2^\omega})$. 
We believe it is an interesting line
of research to further investigate the relationship
between evasion and porosity. In particular, we would
like to know whether some of these cardinals can be shown
to be equal in $ZFC$ (note that $\omega_1 = {\bf e}_X =
unif ({\cal I}_X) < {\bf e}_{2^\omega} = unif ({\cal I}_{
2^\omega}) = \omega_2$ holds in the Mathias real model 
[Br, 3.2.]).
\bigskip
\bigskip
\ce{\sanse References}
\bigskip
\itemitem{[Bl]} {\capit A. Blass,} {\it Cardinal characteristics
and the product of countably many infinite cyclic groups,}
preprint.
\sm
\itemitem{[Br]} {\capit J. Brendle,} {\it Evasion and prediction,}
in preparation.
\sm
\itemitem{[Re]} {\capit M. Repick\'y,} {\it Cardinal invariants
related to porous sets,} to appear in Proceedings of the Bar--Ilan
Conference on Set Theory of the Reals, edited by Haim Judah.
\Smallskip
\no J\"org Brendle \par
\no Department of Mathematics \par
\no Bar--Ilan University \par
\no 52900 Ramat--Gan \par
\no ISRAEL \par
\no {\sanse brendle@bimacs.cs.biu.ac.il} \par
\no (till Sept. 30, 1993)
\bigskip
\no Mathematische Institut \par
\no Universit\"at T\"ubingen \par
\no Auf der Morgenstelle 10 \par
\no 72076 T\"ubingen \par
\no GERMANY \par
\no (after Oct. 1, 1993)

\vfill\eject\end